\documentclass[a4paper,12pt]{article}
\usepackage[T1]{fontenc}
\usepackage[utf8]{inputenc}
\usepackage{lmodern}
\usepackage{graphicx}
\usepackage{amsmath}
\title{\Large\bf Application of Homotopy Perturbation Method  to an Eco-epidemic Model}
\author{P. K. Bera$^{(1)},$ S. Sarwardi$^{(2)}$\footnote{Author to whom all correspondence should be addressed}~ and Md. A. Khan$^{(3)}$\\
$^{(1)}$ Department of Physics, Dumkal College, Basantapur, Dumkal\\ Murshidabad - 742 303, West Bengal, India\\
$^{(2)}$  Department of Mathematics, Aliah University, Salt Lake\\City, Sector-V, Kolkata - 700 091,  West Bengal, India.\\
email: {\bf s.sarwardi@gmail.com}\\
$^{(3)}$  Department of Physics, Aliah University, Salt Lake\\City, Sector-V, Kolkata - 700 091,  West Bengal, India.}
\numberwithin{equation}{section}
\date{}
\begin{document}

\maketitle
%\tableofcontents

\begin{abstract}
In this article, we apply Homotopy Perturbation Method (HPM)
 for solving three coupled non-linear equations
which play an important role in biosystems. To illustrate the capability
and reliability of this method. Numerical example is given which confirms our analytical findings.
\end{abstract}

\section{Introduction}
There exist a wide class of literature dealing with the problem of approximate solutions to
nonlinear equations with various different methodologies, called the perturbation methods. But almost
all perturbation methods are based on small parameters so that the approximate solutions can be expanded
in series of small parameters. Its basic idea is to transform by means of small parameters,
a nonlinear problem of an infinite number of linear subproblems into an infinite number of simpler ones.
 The small parameter determines not only the accuracy of the perturbation approximations but also the
  validity of the perturbation method.

  There exists some analytical approaches, such as the harmonic balance method [1],
  the Krylov-Bogolyubov-Mitropolsky method [2], weighted linearization method [3],
  perturbation procedure for limit cycle analysis [4], modified Lindstedt-Poincare method [5],
  artificial parameter method [6] and so on.

  In science and engineering, there exists many nonlinear problems, which do not contain any
  small parameters, especially those with strong nonlinearity.  He [7, 8] developed the Homotopy
   Perturbation Method (HPM) for solving linear, nonlinear, initial and boundary value problems by merging
     the standard homotopy and the perturbation. The HPM was formulated by taking full advantage of
    the standard homotopy and perturbation methods. In this method the solution is given
    in an infinite series usually converging to an accurate solution.

Inspired and motivated by the ongoing research in the area of bioscience involving mainly ecological and eco-epidemiological systems,  we apply HPM  for solving three coupled nonlinear equations representing a prey-predator model system with disease in prey species only.

      This paper is organized as follows: In Section 2, HPM  has been illustrated. Based on the HPM, the  approximate solutions of three coupled nonlinear equations are obtained in section 3. Finally, we have drawn the conclusion in section 4.

\section{Analysis of the Homotopy Perturbation Method (HPM)}
To illustrate the basic ideas of HPM for solving nonlinear differential
equations, He[7, 8] considered the following nonlinear differential equation:
\begin{equation}
A(u)-f(r)=0, ~\hbox{$r \in \Omega$},\label{eq1}
\end{equation}
with the boundary conditions
\begin{equation}
B\bigg(u,{\partial\over\partial n}\bigg)=0,~ \hbox{$r \in \Gamma$},\label{eq2}
\end{equation}
where $A$ is a general differential operator, $B$ is a boundary operator,
$f(r)$ is known anaclitic function, $\Gamma$ is the boundary of the domain
 $\Omega $ and ${\partial\over \partial n}$ denotes differentiation along
 the normal vector drawn outwards from $\Omega$. The operator $A$ can generally
 be divided into two parts $L$ and $N$, where $L$ is linear and $N$ is nonlinear.
 Therefore, Eq. (\ref{eq1}) can be written as
\begin{equation}
L(u)+N(u)-f(r)=0, r\epsilon \Omega.\label{eq3}
\end{equation}
He [7, 8] constructed a homotopy as follows:
\begin{equation}
H(v,p)=(1-p)\bigl(L(v)-L(u_{0})\bigr)+p\bigl(A(v)-f(r)\bigr)=0\label{eq4}
\end{equation}
or,
\begin{equation}
H_{i}(v,p)=L_{i}(v,p)-L_{i}(v_{0},p)+p\bigl(L_{i}(v_{0},p)+N_{i}(v,p)\bigr)=0.\label{eq5}
\end{equation}

where $v(r,p):\Omega\times[0,1]\to R$. In Eq. (\ref{eq4}), ~$p \in [0,1]$ is an
 embedding parameter and $u_{0}$ is the first approximation that satisfies the boundary
 condition.  The changing process of $p$ from zero to unity is just that of $H(v,p)$
  from $L(v)-L(u_{0})$ to $A(v)-f(r)$. In topology, this is called deformation. The terms
  $L(v)-L(u_{0}) $ and  $A(v)-f(r)$ are called homotopy. According to the homotopy
   perturbation method, the parameter $p$ is used as a small parameter and the solution
   of Eq. (\ref{eq4}) can be expressed as a series in $p$ in the form
\begin{equation}
v=v_{0}+pv_{1}+p^{2}v_{2}+p^{3}v_{3}+\cdots\cdots\cdots\label{eq6}
\end{equation}
 when $p\to 1$, Eq. (\ref{eq4}) corresponds to the original one, Eq. (\ref{eq5})
 the approximate solution  of Eq. (\ref{eq1}) , i.e.
\begin{equation}
u=\lim_{p\to 1}v=v_{0}+v_{1}+v_{2}+v_{3}+\cdots\cdots\cdots\label{eq7}
\end{equation}
The convergence of the series in Eq. (\ref{eq6}) has been discussed by He[7, 8].

\section{Homotopy perturbation method for three coupled system}
We consider three coupled nonlinear equations which describe a prey-predator model, consisting with two prey and a predator species. We also consider an infectious disease, which is transmissible among the prey species only to give the model realism, interested readers are referred to [9-11]. Let us assume $S(t)$ denotes susceptible prey population, $I(t)$ denotes infected prey population, and $P(t)$ denotes predator population at any time. The model under consideration is given by the following system of ordinary nonlinear differential equations
\begin{eqnarray}
{dS\over dt}&=& rS\bigg(1-{S+I\over K}\bigg)-c_{1}SP-\delta SI,\label{eq8}\\
{dI\over dt}&=&\delta SI-c_{2}IP-d_{1}I,\label{eq9}\\
{dP\over dt}&=&e(c_{1}S+c_{2}I)P-d_{2}P,\label{eq10}
\end{eqnarray}
where $ S(0)>0,  I(0)>0, P(0)>0$ and  ${\dot S_{0}}=0,{\dot I_{0}}=0$ and ${\dot P_{0}}=0$.
Here $r$ is the growth rate of the prey population, $c_{1}$ and $c_{2}$ are the searching
efficiency of the predators for the susceptible prey and infected prey respectively,
similarly $ec_{1}$ and $ec_{2}$ are the conversion factors for the susceptible prey and infected
prey respectively consumed by the predators. $K$ is the carrying capacity of the environment for the total (susceptible + infected) prey population, the disease spreads
horizontally with mass action incidence rate $\delta SI$. $d_{1}$ is the mortality rate of infected prey population including disease related death, $d_{2}$ is the mortality rate of
the predator population. All the parameters are non negative.
The predators eat both susceptible and infected prey at different rates, since the
susceptible prey more likely escapes from an attack, thus $c_{1} < c_{2}$. It is to be noted
that the value of the system parameter `e' is a proper fraction for most of the realistic
prey-predator interactions.
For application of HPM, now we write Eqs. (\ref{eq8})-(\ref{eq10}) as
\begin{equation}
H_{i}(S,P,I,p)=L_{i}(S,I,P,p)-L_{i}(S_{0},I_{0},P_{0},p)+p\Bigg(L_{i}(S,I,P,p)+N_{i}(S,I,P,p)\Bigg)\label{eq11}
\end{equation}
where $i=1,2,3$ and
we also consider
\begin{equation}
S=S_{0}+pS_{1}+p^{2}S_{2}+p^{3}S_{3}+\cdots\cdots\cdots\label{eq12}
\end{equation}
\begin{equation}
I=I_{0}+pI_{1}+p^{2}I_{2}+p^{3}I_{3}+\cdots\cdots\cdots\label{eq13}
\end{equation}
\begin{equation}
P=P_{0}+pP_{1}+p^{2}P_{2}+p^{3}P_{3}+\cdots\cdots\cdots\label{eq14}
\end{equation}
 when $p\to 1$,  Eqs. (\ref{eq12})-(\ref{eq14}) become
 the approximate solution  of Eqs. (\ref{eq15})-(\ref{eq17}), i.e.,
\begin{equation}
S_{approx}=\lim_{p\to 1}S=S_{0}+ S_{1}+ S_{2}+ S_{3}+\cdots\cdots\cdots\label{eq15}
\end{equation}
\begin{equation}
I_{approx}=\lim_{p\to 1}I=I_{0}+ I_{1}+ I_{2}+ I_{3}+\cdots\cdots\cdots\label{eq16}
\end{equation}
\begin{equation}
P_{approx}=\lim_{p\to 1}P=P_{0}+ P_{1}+ P_{2}+ P_{3}+\cdots\cdots\cdots\label{eq17}
\end{equation}
Here, boundary conditions are $S_{0}>0, I_{0}>0 ,  P_{0}>0$ and
$\dot S_{0}=\dot I_{0}=\dot P_{0}=0.$

For without perturbation, the Eqs. (\ref{eq8})-(\ref{eq10}) can be written as
\begin{equation}
{dS_{0}(t)\over dt}=rS_{0}\label{eq18}
\end{equation}
\begin{equation}
{dI_{0}\over dt}= -d_{1}I_{0}\label{eq19}
\end{equation}
\begin{equation}
{dP_{0}\over dt}= -d_{2}P_{0}\label{eq20}
\end{equation}
whose solutions are $S_{0}(t)=S(0)e^{rt}, I_{0}(t)=I(0)e^{-d_{1}t}, P_{0}(t)=P(0)e^{-d_{2}t}$.

 With the help of Eq. (\ref{eq11}), one can write the Eqs. (\ref{eq8})-(\ref{eq10}) as follows:
\begin{equation}
{dS\over dt}-rS=p\Bigg( - rS\bigg(  {S+I\over K}\bigg)
-c_{1}SP-\delta SI\Bigg),\label{eq21}
\end{equation}
\begin{equation}
 {dI\over dt}+d_{1}I=p\Bigg( \delta SI-c_{2}IP\Bigg),\label{eq22}
\end{equation}
\begin{equation}
 {dP\over dt}+d_{2}P=p\Bigg( e c_{1}S+ec_{2}I P\Bigg).\label{eq23}
\end{equation}
Substituting the values of $S(t), I(t)$ and $P(t)$ from Eqs. (\ref{eq12})-(\ref{eq14}) and equating the
coefficients of embedding parameter $p$, we get
coefficient of $p^{0}$ as
\begin{eqnarray}
{dS_{0}(t)\over dt}-rS_{0}&=&0,\label{eq24}\\
{dI_{0}\over dt}+d_{1}I_{0}&=&0,\label{eq25}\\
{dP_{0}\over dt} +d_{2}P_{0}&=&0.\label{eq26}
\end{eqnarray}
Equating the coefficient of $p,$ we have
\begin{eqnarray}
 {dS_{1}\over dt}-rS_{1}&=& -{r\over K}S_{0}^{2}-{r\over k}S_{0}I_{0}-c_{1}S_{0}P_{0}-\delta S_{0}I_{0},\label{eq27}\\
{dI_{1}\over dt}+d_{1}I_{1}&=&\delta S_{0}I_{0}-c_{2}I_{0}P_{0},\label{eq28}\\
{dP_{1}\over dt}+d_{2}P_{1}&=&ec_{1}S_{0}P_{0}+ec_{2}I_{0})P_{0},\label{eq29}
\end{eqnarray}
and coefficient of $p^{2}$ as

\begin{eqnarray}
  {dS_{2}\over dt}-rS_{2}&=&-{2r\over K}S_{0}S_{1}-{r\over k}(I_{0}S_{1}+S_{0}I_{1})
 -c_{1}(P_{0}S_{1}+S_{0}P_{1})\nonumber\\&~&-\delta(S_{1}I_{0}+S_{0}I_{1}),\label{eq30}\\
  {dI_{2}\over dt}+d_{1}I_{2}&=& \delta\bigg(S_{1}I_{0}+S_{0}I_{1})-c_{2}(I_{1}P_{0}+I_{0}P_{1}\bigg),\label{eq31}\\
   {dP_{2}\over dt}+d_{2}P_{2}&=& ec_{1}\bigg(S_{1}P_{0}+S_{0}P_{1}\bigg)+ec_{2}\bigg(I_{1}P_{0}+I_{0}P_{1}\bigg),\label{eg32}
\end{eqnarray}
etc.
 Eqs. (\ref{eq8})-(\ref{eq10}) can easily determine the components $S_{k},$ $I_{k}$ and $P_{k}$ and $k\geq 0.$ So, it is possible to calculate more components in the decomposition series to enhance the approximation. Consequently, one can recursively determine every term of the series $\Sigma_{k=0}^{\infty}S_{k}(t),$ $\Sigma_{k=0}^{\infty}I_{k}(t)$ and  $\Sigma_{k=0}^{\infty}P_{k}(t)$ and hence the solutions $S(t)$, $I(t)$ and $P(t)$ is readily obtained in the form of a series like
\begin{eqnarray}S_{approx.}&=&S_{0}e^{rt}+A_{1}e^{2rt}+A_{2}e^{(r-d_{1})t}+A_{3}e^{(r-d_{2})t}
  + A_{4}e^{3rt}+A_{5}e^{(r-d_{1}-d_{2})t}
  \nonumber \\&~& +
   A_{6}e^{(r-2d_{1})t} + A_{7}e^{(r-2d_{2})t}+A_{8}e^{(2r-d_{1})t}+A_{9}e^{(2r-d_{2})t},\label{eg33}\\
 I_{approx.}&=&I_{0}e^{-d_{1}t}+B_{1}e^{(r-d_{1})t}+B_{2}e^{-(d_{1}+d_{2})t}+B_{3}e^{(r-d_{1}-d_{2})t}
  +B_{4}e^{(r-2d_{1})t}\nonumber \\&~&+B_{5}e^{(2r-d_{1})t} +
 B_{6}e^{-(2d_{1}+d_{2})t}+B_{7}e^{-(d_{1}+2d_{2})t},\label{eg34}\\
 P_{approx.}&=&P_{0}e^{-d_{2}t}+C_{1}e^{(r-d_{2})t}+C_{2}e^{-(d_{1}+d_{2})t}+C_{3}e^{(r-d_{1}-d_{2})t}
  +C_{4}e^{(r-2d_{2})t}\nonumber\\&~&+C_{5}e^{(2r-d_{2})t} +
 C_{6}e^{-(2d_{2}+d_{1})t}+C_{7}e^{-(d_{2}+2d_{1})t},\label{eg35}
\end{eqnarray}
where the constants $A_{i}, ~i = 1~ $to$~ 9$ are given by
\begin{eqnarray}A_{1}&=&-{S^{2}(0)\over K},~ A_{2}={1\over d_{1}}\bigg({r\over K}+\delta\bigg)S(0)I(0),~
 A_{3}={c_{1}\over d_{2}}S(0)P(0),\nonumber\\A_{4}&=&{S^{3}(0)\over K^{2}}, ~A_{5}={1\over d_{1}+d_{2}}\Bigg(\bigg({r\over K}+\delta\bigg)\bigg({2c_{2}\over d_{1}}+{c_{1}\over d_{2}}\bigg)-{ec_{2}^{2}\over d_{2}}\Bigg)S(0)I(0)P(0),\nonumber\\A_{6}&=&{rc_{1}\over 2Kd_{1}d_{2}}I(0)S(0)P(0)+{\delta\over 2d_{1}^{2}}\bigg({r\over K}+\delta\bigg)I^{2}(0)S(0),\nonumber\\A_{7}&=&{c_{1}^{2}\over d_{2}^{2}}S(0)P^{2}(0), ~A_{8}=-{1\over r-d_{1}}\bigg({\delta\over r}+{1\over K}\bigg)\bigg({2r^{2}\over Kd_{1}}-{r\over K}+\delta\bigg)S^{2}(0)I(0),\nonumber\\A_{9}&=&-{1\over r-d_{2}} \bigg({2rc_{1}\over Kd_{2}}-{ec_{1}^{2}\over r}-{c_{1}\over K}\bigg)S^{2}(0)P(0),\nonumber
\end{eqnarray}
the constants $B_{i},~i = 1~ $ to $~ 7$ are given by
\begin{eqnarray}B_{1}&=&{\delta\over r}S(0)I(0),~B_{2}={c_{2}\over d_{1}}I(0)P(0),~B_{3}={1\over r-d_{2}}\Bigg(\delta\bigg({c_{1}\over d_{2}}+{c_{2}\over d_{1}}\bigg)\nonumber\\&~&-c_{2}\bigg({ec_{1}+\delta\over r}\bigg)\Bigg) S(0)I(0)P(0),\nonumber\\B_{4}&=&{\delta c_{1}\over d_{2}(r-d_{1})}S(0)I(0)P(0),B_{5}={\delta\over 2r}\Bigg({\delta\over r}-{1\over K}\Bigg)S^{2}(0)I(0),\nonumber\\B_{6}&=&-{ec_{2}^{2}\over d_{2}(d_{1}+d_{2})}I^{2}(0)P(0),B_{7}={c_{2}^{2}\over 2d_{2}d_{1}}P^{2}(0)I(0),\nonumber
\end{eqnarray}
and the constants $C_{i},~i = 1~ $to$~ 7$ are given by
\begin{eqnarray}C_{1}&=&{ec_{1}\over r}S(0)P(0),~ C_{2}=-{ec_{2}\over d_{2}}I(0)P(0),\nonumber\\C_{3}&=&{1\over r-d_{1}}\Bigg(\bigg({e^{2}c_{1}c_{2}\over r}+{ec_{2}\delta\over r}-{e^{2}c_{1}c_{2}\over d_{2}}\bigg)S(0)I(0)P(0)+{ec_{1}^{2}\over d_{2}}S(0)P^{2}(0)\Bigg),\nonumber\\C_{4}&=&{ec_{1}^{2}\over d_{2}(r-d_{2})}S(0)P^{2}(0),
 ~C_{5}={ec_{1}\over 2r}\Bigg({ec_{1}\over r}-{1\over K}\Bigg)S^{2}(0)P(0),\nonumber\\ C_{6}&=&-{ec_{2}^{2}\over d_{1}(d_{1}+d_{2})}P^{2}(0)I(0),
 ~C_{7}={e^{2}c_{2}^{2}\over 2d_{2}d_{1}}I^{2}(0)P(0).\nonumber
\end{eqnarray}

Putting the different values of parameters $r,k,c_{1},c_{2},\delta,e,d_{1},d_{2},$ we obtained the approximate solutions of Eqs. (\ref{eq8})-(\ref{eq10}). To explain our results, we have   drawn   figures of $S(t), I(t) $ and $P(t)$ versus time for different parameters and initial condition values. From these figures, we see that the exact numerical results and the results obtained using HPM are nearly same.

\section{Conclusion}
This system of three coupled differential equations (3.1)--(3.3) plays an important role in biosystems. The basic goal of
 this paper is to study this model using HPM. The goal has been achieved by deriving solutions using few iterations only. The qualitative results of the present studies have been compared with the results obtained by  numerical computation using  $r=0.1,k=0.3,c_{1}=0.1,c_{2}=0.2,\delta=0.1,e=0.1,d_{1}=0.2,d_{2}=0.2$ as evident from the Fig. 1  and reveal that HPM is very effective and convenient for solving non-linear differential equations. Hope that with the help of  these solutions, one can study qualitative and quantitative behaviors of  realistic prey-predator interactions. The  HPM introduces a significant improvement in this field. This makes the proposed scheme more powerful and gives a wider applicability.

 \begin{figure}[ht]\label{f1}
\begin{center}
\includegraphics[height = 8cm, width=12cm]{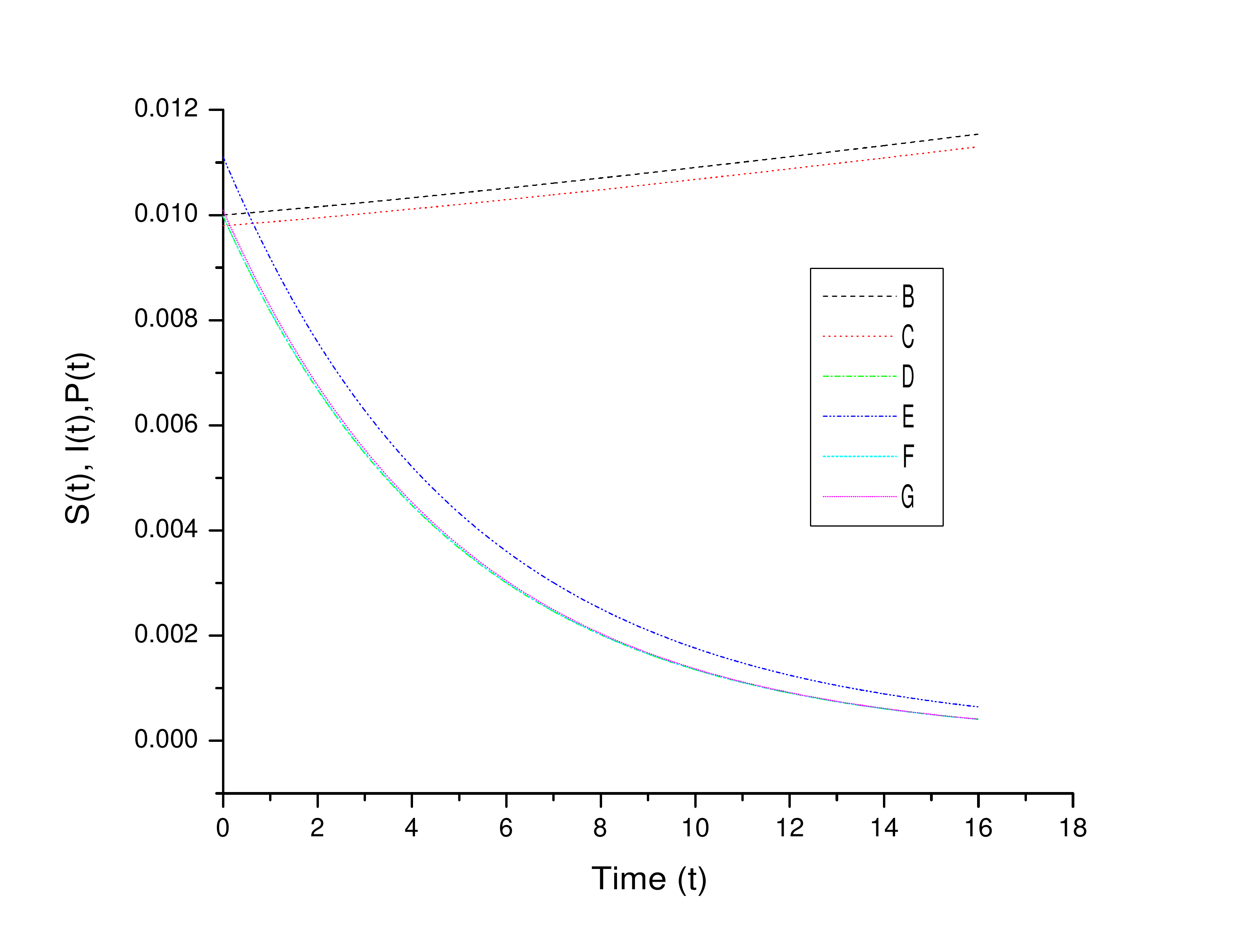}
\caption{\textrm{\small The $S(t),I(t), P(t)$  versus time for  the different values of parameters $r=0.1,k=0.3,c_{1}=0.1,c_{2}=0.2,\delta=0.1,e=0.1,d_{1}=0.2,d_{2}=0.2$ and initial conditions $S_{0}(0)=0.01, I_{0}(0)=0.01$ and $P_{0}(0)=0.01$. B-, D-, F-line represent the numerical solutions $S(t), I(t)$ and $P(t)$ respectively and  C-, E-, G-line represent the approximate solutions of $S(t), I(t)$ and $P(t)$ respectively which have been  obtained using HPM.}}
\end{center}
\end{figure}

\vskip7.0cm

{\bf Acknowledgement:} This work is  supported by University Grants Commission(UGC), Government of  India (Project No.PSW-63/12-13(ERO)). Dr. S. Sarwardi is thankful to the Department of Mathematics, Aliah University for extending opportunities to perform the present work.

\end{document}